\documentclass[12pt]{article}
\usepackage[T2A]{fontenc}
\usepackage{amssymb}
\usepackage{amsmath}
\usepackage{euscript}

\begin{document}

\centerline{\large\bf On the Hyperbolicity Properties of Inertial
Manifolds} \centerline{\large\bf of Reaction--Diffusion Equations}
\medskip
\centerline{\bf A. V. Romanov}
\medskip
\centerline{\it National Research University Higher School of
Economics, Moscow, Russia} \centerline{\it e-mail:
av.romanov@hse.ru}
\bigskip

\noindent\textit{Abstract}. For 3D reaction--diffusion equations,
we study the problem of existence or nonexistence of an inertial
manifold that is normally hyperbolic or absolutely normally
hyperbolic. We present a system of two coupled equations with a
cubic nonlinearity which does not admit a normally hyperbolic
inertial manifold. An example separating the classes of such
equations admitting an inertial manifold and a normally hyperbolic
inertial manifold is constructed. Similar questions concerning
absolutely normally hyperbolic inertial manifolds are discussed.

\bigskip



\bigskip
\medskip
\centerline{\bf  0. Introduction}
\bigskip

The existence of a smooth inertial manifold ${\EuScript M}$ for
the dissipative parabolic equation in the infinite-dimensional
Hilbert space implies [14,18,19] that its final dynamics (as
$t\rightarrow+\infty$) is controlled by finitely many parameters.
The additional property of normal hyperbolicity of the inertial
manifold ${\EuScript M}$ guarantees the structural stability of
this manifold. The stronger property of absolute normal
hyperbolicity means one and the same hyperbolicity parameters for
the entire ${\EuScript M}$. So far, the existence of an inertial
$C^{1}$-manifold has been established for a rather narrow class of
semilinear parabolic equations, while known examples of its
nonexistence [2,15,16] seem to be somewhat artificial and are not
related to problems of mathematical physics.

The present paper deals with necessary conditions for the
existence of the above-mentioned two types of inertial manifolds
of scalar and vector reaction--diffusion equations. For the 3D
chemical kinetics equations with a cubic nonlinearity, we strive
for constructing examples separating the classes of problems
admitting an \textit{inertial manifold}, a \textit{normally
hyperbolic inertial manifold}, and an \textit{absolutely normally
hyperbolic inertial manifold}. An example separating the first two
possibilities is obtained for two-component systems. Namely, in
Proposition 3.5 we construct an (uncoupled) system of such
equations that has an inertial manifold but does not admit a
normally hyperbolic inertial manifold. In particular, this system
provides an example of an inertial manifold that is not normally
hyperbolic. On the other hand, we present a system of two
\textit{coupled} reaction--diffusion equations of this type that
do not admit a normally hyperbolic inertial manifold in the
natural state space (Proposition 3.4). An example of a scalar 3D
equation with a cubic nonlinearity without an absolutely normally
hyperbolic inertial manifold is constructed. Note that the order
of the polynomial nonlinearity in the chemical kinetics equations
corresponds to the \textit{reaction order}, which usually does not
exceed~$3$. We also discuss how close the well-known sufficient
conditions (the \textit{spectral jump condition} and the
\textit{spatial averaging principle}) for the existence of
strongly and weakly normally hyperbolic inertial manifolds are to
being necessary.

The paper is organized as follows. Section~1 contains elementary
information about abstract semilinear parabolic equations. The
necessary and sufficient conditions, known so far, for the
existence of a smooth inertial manifold are stated in Section~2.
The main results on the existence and nonexistence of various
inertial manifolds for the reaction--diffusion equations are
presented in Sections~3--4. Section~5 discusses conjectures on the
relationship between spectral properties of the linear part of the
equation and the existence or nonexistence of various types of
inertial manifolds.

The results of the paper were presented by the author at the
Conference-School ``Infinite-dimensional dynamics, dissipative
systems, and attractors'' held at the Lobachevsky State University
of Nizhny Novgorod on July 13--17, 2015.

\bigskip
\centerline{\bf 1. Preliminaries}
\bigskip

A semilinear parabolic equation in a real separable
infinite-dimensional Hilbert space $(X,\left\| \, \cdot \,
\right\|)$ has the form
$$
  \partial _{t} u \, = -Au+F(u).\eqno (1.1)
$$
Here we assume that

(i) $A:{\EuScript D}(A)\rightarrow X$ is a linear positive definite
self-adjoint operator with compact inverse $A^{-1}$.

(ii) $F\in C^{1}(X^{\theta},X)$ is a nonlinear function with domain
$X^{\theta}={\EuScript D}(A^{\theta})$, $0\leq\theta<1$,
$\|u\|_\theta=\|A^{\theta}u\|$, such that
$$
   \|F(u_{1})-F(u_{2})\|\leq L(r)\|u_{1}-u_{2}\|_\theta \eqno (1.2)
$$
on the balls ${\EuScript B}_{r}=\{u\in X^{\theta}:
\|u\|_\theta<r\}$.

(iii) There exists a dissipative phase semiflow
$\{\Phi_{t}\}_{t\geq0}$ on $X^{\theta}$.

We refer to the number $\theta$ as the \textit{nonlinearity
exponent} of Eq.~(1.1) and set $X^{0}=X$. The space $X$ will be
called the \textit{main} space. Dissipativity is understood as the
existence of an absorbing ball ${\EuScript B}_{r}\subset
X^{\theta}$ (see~[14,~19]). Under these conditions~[4], the phase
semiflow proves to be smooth, and the evolution operators
$\Phi_{t}:X^{\theta}\rightarrow X^{\theta}$, $t>0$, are compact.
The \textit{parabolic smoothing} property guarantees the inclusion
$\Phi_{t}X^{\theta}\subset X^{1}={\EuScript D}(A)$ for $t>0$.

The global attractor ${\EuScript A}$ is defined as the union of all
complete bounded trajectories of the equation; in our case, it is a
compact subset of $X^{\theta}$. An \textit{inertial manifold} of
Eq.~(1.1) is a \textit{smooth} ($C^{1}$) finite-dimensional
positively invariant surface ${\EuScript M}\subset X^{\theta}$
containing the attractor ${\EuScript A}$ and attracting all
trajectories $u(t)$ with exponential tracking as
$t\rightarrow+\infty$. An inertial manifold usually has a Cartesian
structure and is diffeomorphic to a ball in ${\mathbb R}^{n}$. The
restriction of (1.1) to ${\EuScript M}$ gives an inertial form (an
ordinary differential equation in ${\mathbb R}^{n}$, $n=\dim
{\EuScript M}$), which completely reproduces the final dynamics of
the original equation. There is a vast literature dealing with the
theory of inertial manifolds (see [14,18--20] and references
therein); moreover, one often considers \textit{Lipschitz}
(nonsmooth) inertial manifolds.

\bigskip
\centerline{\bf 2. Inertial Manifold: Existence Conditions}
\bigskip

The dissipativity of the evolution system (1.1) permits one to
change the function $F(u)$ outside ${\EuScript B}_{r}$ with the
preservation of $C^{1}$-regularity in such a way that the new
function $\tilde{F}(u)$ is identically zero outside the ball
${\EuScript B}_{r+1}$. This ``truncation'' procedure (e.g., see
[19]) permits one to proceed to the equation
$$
 u_{t} =-Au+\tilde{F}(u), \eqno (2.1)
$$
which inherits the final dynamics of the original problem. One has
$L(r)\equiv L$ in the estimate (1.2) for $\tilde{F}(u)$. It is well
known [1, 17, 18] that the existence of a smooth $n$-dimensional
inertial manifold ${\EuScript M}\subset X^{\theta}$ of Eq.~(2.1) in
the phase space $X^{\theta}$ is guaranteed by the spectral jump
condition $\mu_{n+1}-\mu_{n}>cL(\mu_{n+1}^{\theta}
+\mu_{n}^{\theta})$, where $0<\mu _{1}\leq\mu_{2}\dotsm$ are the
eigenvalues of the operator~$A$ arranged in nondescending order
(counting multiplicities) and $c>0$ is an absolute constant. The
manifold ${\EuScript M}$ also proves to be an inertial manifold of
the original parabolic equation. Thus, the \textit{spectrum
sparseness} condition
$$
 \mathop{{\rm sup}}\limits_{n\ge 1} \, \, \frac{\mu _{n+1} -\mu _{n} }{\mu _{n+1}^{\theta }
 +\mu _{n}^{\theta } } \;  =\,  \infty \eqno (2.2)
$$
is sufficient for the existence of an inertial $C^{1}$-manifold
${\EuScript M}\subset X^{\theta}$ of the dissipative equation (1.1)
with given linear part $-A$ for an arbitrary nonlinear function
$F:X^{\theta} \rightarrow X$ with properties (ii).

Now consider the scalar reaction--diffusion equation
$$
  \partial _{t} u \, = \, \nu\Delta u
  +  f(x,u),\;\,\nu>0, \eqno (2.3)
$$
in a bounded Lipschitz domain $\Omega \subset {\mathbb R}^{m}$
with one of the standard boundary conditions $(\rm{D})$,
$(\rm{N})$, or $(\rm{P})$ and with a sufficiently smooth function
$f:\bar{\Omega}\times\mathbb{R}\rightarrow \mathbb{R}$ satisfying
the \textit{sign condition} $v\cdot f(x,v)<0$ for $x\in \Omega$
and $|v|\geq r>0$. In this case, there exists a dissipative phase
semiflow on $X=L^{2}(\Omega)$ [19, Chapter~3]. Let us extend
$f(\cdot,v)$ from $[-r,r]$ with the preservation of smoothness to
a Lipschitz function $\tilde{f}(\cdot,v)$ vanishing for $|v|\geq
r+1$. By the maximum principle, the partial differential
equation~(2.3) with $f(x,u)$ replaced by $\tilde{f}(x,u)$ inherits
the limit modes of the original problem and admits the
interpretation (1.1) with nonlinearity exponent $\theta =0$ and
with $X^{1}\subset H^{2}(\Omega)$. To this end, one should set
$Au=u-\nu\Delta u$ and $F(u) =u+\tilde{f}(x,u)$.

For $m\leq3$, the well-known difficulties [4, p.~11] concerning
the smoothness of the Nemytskii operator in $L^{2}(\Omega)$ force
one to use the weakened version ([11, pp.~813, 836]; see also
[15]) of the definition of Fr\'echet derivative of the nonlinear
function $u\rightarrow F(u)$, where one requires that  $u,h \in
X^{1}$ in the analysis of the increment $F(u+h)-F(u)$. This
approach (\textit{generalized Fr\'echet derivative}), which uses
the parabolic smoothing property, was generalized in [6,
Section~7]. The phase semiflow of Eq.~(2.3) in $X=L^{2}(\Omega)$
is differentiable in the same sense.

If the spectrum is
$\sigma(-\Delta)=\{0\leq\lambda_{1}\leq\lambda_{2}\leq\dotsm\}$,
then condition (2.2) is reduced to the relation
$$
 \mathop{{\rm sup}}\limits_{n\ge 1} \,
 (\lambda _{n+1}-\lambda _{n} )=\infty \, ,   \eqno (2.4)
$$
which seems to be rather restrictive in view of the Weyl
asymptotics $\lambda _{n} \sim c\, n^{2/m}$. We point out that
(2.4) holds for $m=1$ as well as for some domains $\Omega \subset
{\mathbb R}^{2} $. These domains include rectangles with rational
squared side ratio [11], but in general the description of planar
domains for which $\sigma (-\Delta )$ is sparse remains a mystery.
Already for $m=3$, one has $\lambda _{n} \sim cn^{2/3} $, and
condition (2.4) seems to be exotic.

In this connection, the following property of the Laplace operator
in a domain $\Omega \subset {\mathbb R}^{m}$, $m\leq3$, was stated
in~[10, 11], which was referred there to as the \textit{principle
of spatial averaging}. Set
$$
  (B_{h}u)(x)=h(x)u(x),\;\;\;\;\bar{h}=({\rm vol} \,
\Omega)^{-1}\int \limits _ \Omega h(x)dx
$$
for $h\in H^{2}(\Omega) \subset L^{\infty}(\Omega)$ and $u\in
L^{2}(\Omega)$. Let $P_{\lambda}$ be the spectral projection of the
self-adjoint operator $-\Delta$ corresponding to the part of the
spectrum in $[0,\lambda]$, and let $I=\rm id$.

DEFINITION 2.1. \textit{The Laplace operator $\Delta$ with a given
standard boundary condition in a bounded Lipschitz domain $\Omega
\subset {\mathbb R}^{m}$, $m\leq3$, satisfies the principle of
spatial averaging if there exists a $\rho >0$ such that for any
$\varepsilon >0$ and $k>0$ there exists an arbitrarily large
$\lambda>k$ such that $\lambda\in [\lambda_{n},\lambda_{n+1})$,
$\lambda_{n+1}-\lambda_{n}\geq\rho$, and
$$
  \|(P_{\lambda+k}-P_{\lambda-k})(B_{h}-\bar{h}I)
  (P_{\lambda+k}-P_{\lambda-k})\|_{\rm op}
  \; \leq \; \varepsilon \|h\|_{H^{2}}\;\:\:\:\forall \;h\in
   H^{2}(\Omega), \eqno (2.5)
$$
where $\|\cdot\|_{\rm op}$ is the norm on ${\rm End} \,
L^{2}(\Omega)$.}

Essentially, one speaks of an arbitrarily good approximation, for
any $h\in H^{2}(\Omega)$, to the Schr\"odinger operator
$\Delta+h(x)I$ by a shifted Laplace operator $\Delta+\bar{h}I$ in
an arbitrarily wide range of eigenmodes of the Laplace operator.
Here one assumes that
$$
{\rm lim}\mathop{{\rm
sup}}\limits_{n\rightarrow\infty}(\lambda_{n+1}-\lambda_{n})>0,
$$
which is always the case for $m\leq2$. This principle follows from
the sparseness of the spectrum  (but not vice versa!) and ensures
[11, p.~846] the existence of a smooth inertial manifold of
Eq.~(2.3) with $f\in C^{3}$. In particular, the principle of
spatial averaging holds for an \textit{arbitrary} rectangle
$\Omega_{2} \subset {\mathbb R}^{2} $ and for a cube $\Omega_{3}
\subset {\mathbb R}^{3} $ [11], although condition (2.4) is not
guaranteed for the former and is violated for the latter. In~[8],
the existence of a (Lipschitz) inertial manifold of Eq.~(2.3) is
derived from less restrictive conditions: the number $\lambda>k$
may depend on bounded sets ${\EuScript B}\subset H^{2}(\Omega)$,
and (2.5) is replaced by the estimate
$$
 \|(P_{\lambda+k}-P_{\lambda-k})(B_{h}-\bar{h}I)(P_{\lambda+k}
 -P_{\lambda-k})\|_{\rm
  op} \; \leq \; \varepsilon \;\:\:\:\forall \;h\in {\EuScript B}.
$$
In the framework of this approach, the existence of an inertial
manifold was proved for Eq.~(2.2) in some 2D and 3D
polyhedra~[8,~9]. The principle of spatial averaging has only been
proved to hold in some model cases, and unfortunately, this
principle practically does not apply to systems of
reaction--diffusion equations, because in this case the operator
corresponding to the componentwise multiplier is the operator of
multiplication by a matrix of numbers that is diagonal \textit{but
not scalar}.

Recently, Zelik [20] suggested an abstract form of the principle of
spatial averaging, which generalizes the constructions in [8--11]
and ensures the existence of a smooth inertial manifold of
Eq.~(1.1). This approach was further developed in~[6, 7]. The
corresponding technique permitted establishing the existence of an
inertial manifold ${\EuScript M}\in C^{1+\varepsilon}$ for the
Cahn--Hilliard equation [6] and of an inertial manifold ${\EuScript
M}\in {\rm Lip}$ for the modified Leray $\alpha$-model of the
Navier--Stokes equations on the three-dimensional torus [7].

So far, little is known about the cases of nonexistence of an
inertial manifold for parabolic problems. A system of two coupled
one-dimensional parabolic pseudodifferential equations that does
not admit a smooth inertial manifold was constructed in~[15]. A
general construction of abstract equations (1.1) with nonlinearity
exponent $\theta =0$ and without a smooth inertial manifold is
described in~[2]. A more natural  story is considered in [16],
where an integro-differential parabolic equation with nonlocal
diffusion on the circle is presented which does not have an
inertial manifold in the chosen state space.

All these examples are based on the following argument. Since the
phase semiflow is dissipative and compact, it follows that the
stationary point set $E=\{u\in X^{1}:F(u)-Au=0\}$ of Eq.~(1.1) is
nonempty. Since $E\subset {\EuScript A}$, we see that $E$ is
contained in the inertial manifold, provided that the latter
exists. Since the operator $A^{-1}$ is compact and, by [4,
Chapter~1], the linear operator $-S_{u}=A-F'(u)$ on $X$ is
sectorial, it follows that the spectrum $\sigma (S_{u})$, $u\in
E$, consists of eigenvalues $\lambda$ of finite multiplicity, and
the number $l(u)$ (counting multiplicities) of positive $\lambda$
in $\sigma (S_{u})$ is finite. Let $E_{-}=\{u\in E: \sigma
(S_{u})\cap (-\infty,0]=\phi\}$.

Now we can state a necessary condition for the existence of an
inertial manifold as follows.

LEMMA 2.2 ([15]). \textit{If Eq.~$(1.1)$ admits a smooth inertial
manifold ${\EuScript M}\subset X^{\theta} $}, \textit{then the
number $l(u_{0} )- l(u_{1} )$ is even for any $u_{0},u_{1} \in
E_{-}$.}

To apply the lemma, one usually constructs a nonlinearity $F$ such
that Eq.~(1.1) has stationary solutions $u_{0},u_{1}\in E_{-}$ with
$l(u_{0})=0$ and $l(u_{1})=1$.

\bigskip
\centerline{\bf 3.  Normally Hyperbolic Inertial Manifolds}
\bigskip

Unfortunately, so far there are no examples physically more
meaningful than those given above of parabolic equations without
inertial manifolds. At the same time, such examples were obtained
in~[12,~15] for the case in which one speaks of inertial manifolds
with additional hyperbolicity properties.

DEFINITION 3.1. A \textit{smooth inertial manifold ${\EuScript
M}\subset X^{\theta}$ of Eq.}~(1.1) \textit{is said to be normally
hyperbolic if, for some vector bundle ${\EuScript T}_{{\EuScript
M}}X^{\theta}={\EuScript T}{\EuScript M} \oplus {\EuScript N}$
invariant with respect to the linearization $\{\Phi '_{t}\}$ of
the semiflow $\{\Phi_{t}\}_{t\geq0}$, where ${\EuScript
T}{\EuScript M}$ is the tangent bundle of ${\EuScript M}$, one has
the estimates
$$\left\| \, \Phi '_{t} (u)h \, \right\| _{\theta}  \, \ge
\, M^{-1} e^{-\gamma _{1} t} \left\| \, h \, \right\| _{\theta} \;
\ (h \in {\EuScript T}_{u}{\EuScript M}),$$
$$\left\| \, \Phi '_{t} (u)h \, \right\| _{\theta}  \, \le
\, M e^{-\gamma _{2} t} \left\| \, h \, \right\| _{\theta} \; \ (h
\in {\EuScript N}_{u}) \eqno (3.1)$$
with constants $M>0$ and $0 <
\gamma_{1} < \gamma_{2}$ depending on ${\EuScript M}$ and $u\in
{\EuScript M}$. If these constants are independent of $u\in
{\EuScript M}$, then the manifold is said to be absolutely
normally hyperbolic}.

We point out that the normally hyperbolic invariant manifolds of
finite- and infinite-dimensional dynamical systems are structurally
stable~[5,~13].

The methods in~[6] permit one to establish that the validity of the
abstract version of the principle of spatial averaging [20] implies
the existence of a normally hyperbolic inertial manifold in the
state space of the parabolic problem~(1.1). For the
reaction--diffusion equations (2.3), as similar claim was announced
as early as in~[10; 11, p.~830].

The known necessary conditions for the existence of an inertial
manifold ${\EuScript M}\subset X^{\theta}$ with hyperbolicity
properties amount to analyzing the spectrum of the linearization of
the vector field $F(u)-Au$ of Eq.~(1.1) on the stationary point set
$E\subset X^{1}$. For $\gamma\in \mathbb{R}$ and $u\in E$, let
$Y(u,\gamma )$ be the finite-dimensional invariant subspace of the
operator $S_{u}=F'(u)-A$ corresponding to the part of the spectrum
$\sigma (S_{u})$ with ${\rm Re}\, \lambda  \geq \gamma $.

Lemma 3.2 ([12, 15]). \textit{If the inertial manifold ${\EuScript
M}\subset X^{\theta}$ of Eq.}~(1.1) \textit{is normally hyperbolic,
then
$$
 \forall \,u\in E,\; \, \exists \,\gamma=\gamma(u;{\EuScript
M})<0:\;\,{\rm dim}\,Y(u,\gamma)={\rm dim}\,{\EuScript M}.
$$
In the case of absolutely normal hyperbolicity of ${\EuScript
M}\subset X^{\theta}$, one has $ \gamma=\gamma({\EuScript M})$.}

Here $\gamma=-(\gamma_{1}+\gamma_{2})/2$, where
$0<\gamma_{1}<\gamma_{2}$ are the numbers in Definition~3.1. For
$u\in E$, the invariant subspaces ${\EuScript T}_{u}{\EuScript M}$
and ${\EuScript N}_{u}$ of the operator $S_{u}$ correspond to the
parts of the spectrum $\sigma(S_{u})$ with ${\rm Re}\,\lambda\geq
-\gamma_{1}$ and ${\rm Re}\,\lambda\leq -\gamma_{2}$, respectively;
moreover, $\Phi '_{t} (u)={\rm exp}\; (-tS_{u})$, $t>0$.

The lemma was used to obtain the well-known example [12,
Theorem~2.5] of Eq.~(2.3) in the cube $\Omega =(0,\pi )^{4} $ with
the Neumann condition on $\partial \Omega$ and with a real-analytic
function $f(x,u)$ (polynomial in $u$) for which there does not
exist a normally hyperbolic inertial manifold ${\EuScript M}\subset
L^{2} (\Omega )$. However, the function $f$ was not constructed in
closed form in this example. Furthermore, it would be of interest
to obtain similar examples for 3D reaction--diffusion equations
with a homogeneous polynomial nonlinearity~$f(u)$. Moreover, from
the viewpoint of applications to chemical kinetics, the degrees of
the polynomials should not exceed~$3$.

Consider the two-component system
$$
\partial _{t} u_{1} \, =\, \Delta u_{1} \, +
f_{1}(u_{1},u_{2}),\;\;\; \partial _{t} u_{2} \, =\, \Delta u_{2}
\, + f_{2}(u_{1},u_{2})\eqno (3.2)
$$
in the cube $\Omega =(0,\pi )^{3} $ with the Neumann condition
$(\rm{N})$ on $\partial \Omega$ and with a $C^{3}$-function
$f=(f_{1},f_{2})$, ${\mathbb R}^{2}\stackrel{f}{\longrightarrow}
{\mathbb R}^{2}$. Then, just as above, system (3.2) can be reduced
to the abstract dissipative problem (1.1) with
$X=L^{2}(\Omega;{\mathbb R}^{2})$ and with the nonlinearity
exponent $\theta=0 $ under the assumption that there exists an
\textit{invariant region} [19, Chapter~3] for the ordinary
differential equation $v_{t}=f(v),\, v\in {\mathbb R}^{2}$. Here
the smoothness of the operator $u\rightarrow f(u)$, $u\in X$, is
understood in the sense of the weakened Fr\'echet derivative.

For a fixed point $p\in {\mathbb R}^{2}$ of the vector field $f$,
we set $\delta(p)=|{\rm Re} \, (\xi_{1}-\xi_{2})|$, where
$\xi_{1}$ and $\xi_{2}$ are the eigenvalues of the Jacobian matrix
$f'(p)$. Note that $\delta(p)=0$ in the case of multiple or
complex eigenvalues of the matrix~$f'(p)$.

LEMMA 3.3 ([15]). \textit{The dissipative system} $(3.2)$
\textit{does not have a normally hyperbolic inertial manifold in
the state space~$X$ if the vector field~$f$ has four fixed points
$p_{i} \in {\mathbb R}^{2}$ such that $\delta(p_{i})=i$ for
$i=0,1,2,3$}.

The proof uses the necessary condition given by Lemma~3.2. The
existence of a smooth vector field $f$ with the desired properties
on ${\mathbb R}^{2}$ is obvious. Our aim is to construct a
third-order polynomial field of this kind. Set
$$
  f_{1}(v_{1},v_{2})=kv_{1}(1-av_{1}^{2}+v_{2}^{2}), \; \;
  f_{2}(v_{1},v_{2})=kv_{2}(1-bv_{2}^{2}-v_{1}^{2}) \eqno (3.3)
$$
with some constants $k,a,b>0$.

We have $v\cdot f(v)\leq 0$ for $|v|^{2}\geq r_{0}^{2}=2/{\rm
min}(a,b)$ and dissipativity of the system (3.2) with the vector
field $(3.3)$ is ensured by the positive invariance of the disks
$|v|\leq r$ with $r\geq r_{0}$ for the ordinary differential
equation $v_{t}=f(v)$, $v\in {\mathbb R}^{2}$. Furthermore, the
condition $b^{-1} \leq c^{2} \leq a-1$ implies the positive
invariance of the region $D_{c}=\{v\in {\mathbb R}^{2}: 0 \leq
v_{1} \leq 1,\, 0\leq v_{2} \leq c\}$ for the equation
$v_{t}=f(v)$, which in its turn imply [19] the preservation of
this region for the components $u_{1},u_{2}$ in the system (3.2).

Proposition 3.4. \textit{There exist positive $k,a$ and $b$ with
$a\geq 1+b^{-1}\geq b$, such that the dissipative coupled system
$(3.2)$ with the vector field $(3.3)$ and with $\Omega =(0,\pi
)^{3}$ does not have a normally hyperbolic inertial manifold
${\EuScript M}\subset X$}.

PROOF. Assuming that $a>1$, let us single out four fixed points
$$p_{0}=(0,0), \;\,p_{1}=(\frac{1} {\sqrt{a}},0), \;\, p_{2}= (\sqrt{\frac{b+1}{ab+1}},
\sqrt{\frac{a-1}{ab+1}}), \;\, p_{3}=(0,\frac{1} {\sqrt{b}})$$ of
the vector field $f$ on ${\mathbb R}^{2}$. Here
$$
f'(v)=k\left(\begin{array}{cc} {1-3av_{1}^{2}+v_{2}^{2}} &
{2v_{1}v_{2}} \\ {-2v_{1}v_{2}} & {1-v_{1}^{2}-3bv_{2}^{2}}
\end{array}\right)
$$
for $v\in {\mathbb R}^{2}$ and
$$
f'(p_{0})=\left(\begin{array}{cc} {k} & {0} \\ {0}
& {k}
\end{array}\right), \; \; \; \;
f'(p_{1})=\left(\begin{array}{cc} {-2k} & {0} \\ {0}
& {k-k/a}
\end{array}\right),
$$
$$
f'(p_{2})=k\left(\begin{array}{cc} {\frac{-2ab-2a}{ab+1}} &
{{\frac{2((a-1)(b+1))^{1/2}}{ab+1}}} \\
{{\frac{-2((a-1)(b+1))^{1/2}}{ab+1}}} & {\frac{-2ab+2b}{ab+1}}
\end{array}\right), \; \; \; f'(p_{3})=\left(\begin{array}{cc}
{k+k/b} & {0} \\ {0} & {-2k}
\end{array}\right).
$$
\medskip
Set $\delta_{i}=\delta(p_{i})$, $0\leq i\leq 3$. We have
$\delta_{0}=0$, $\delta_{1}=k(3-a^{-1})$,
$\delta_{3}=k(3+b^{-1})$, and
$$
\delta_{2}^{2}= \frac{4k^{2}(a+b)^{2}}{(ab+1)^{2}}-16k^{2}
\frac{(a-1)(b+1)}{(ab+1)^{2}},
$$
$$\delta_{2}=\frac{2k}{ab+1}\cdot |a-b-2|.$$
Set $k=a/(3a-1)$ and $b=a/(6a-3)$; then $\delta_{1}=1$ and
$\delta_{3}=3$. The function $\varphi:a\rightarrow\delta_{2}$ is
continuous on $(1,\infty)$, and, since $k(\infty)=1/3,\,
b(\infty)=1/6$, we have $\varphi(7)<2,\, \varphi(\infty)=4$. Thus,
there exists $a=a^{\ast}>7$ such that $\varphi(a)=2$. It is easy
to verify that: 1) $a\geq 1+b^{-1}\geq b$ and $r_{0}^{2}=2/b<12$;
2) $p_{i} \in D_{c}$ and $|p_{i}|\leq \sqrt{7} $ for $c=\sqrt{6}$
and $0\leq i \leq3$. Since $\delta(p_{i})=i,\, 0 \leq i \leq 3$,
the proposition follows from Lemma 3.3. $\Box$

Now consider the vector field
$$
 f_{1}(v_{1},v_{2})=v_{1}(a-v_{1})(v_{1}-b), \; \;
 f_{2}(v_{1},v_{2})=v_{2}(c-v_{2})(v_{2}-d) \eqno (3.4)
$$
with $a=2$, $b=\sqrt{3}$, $c=\sqrt{6}$, and $d=\sqrt{2}$. The
dissipativity and the preservation of the positivity of solutions
of the corresponding problem~(3.2) is guaranteed by the sign
condition with respect to each component and by the positive
invariance of the quadrant $v_{1}\geq 0,v_{2}\geq 0$ with respect
to the ordinary differential equation $v_{t}=f(v)$ in ${\mathbb
R}^{2}$.

PROPOSITION 3.5. \textit{The dissipative uncoupled system $(3.2)$
with the vector field $(3.4)$ and with $\Omega =(0,\pi )^{3}$
admits an inertial manifold ${\EuScript M}\subset X$ but does not
have a normally hyperbolic inertial manifold in~$X$}.

PROOF. Each of the scalar equations in~(3.2) admits an inertial
manifold ${\EuScript M_{j}}\subset L^{2}(\Omega)$, $j=1,2$ [11],
and hence ${\EuScript M}={\EuScript M_{1}}\times {\EuScript M_{2}}$
is an inertial manifold of the two-component system in $X$. At the
stationary points
$$
 p_{0}=(0,0), \; \, p_{1}=(b,d), \; \, p_{2}=(a,c), \;
 \, p_{3}=(b,c),
$$
the Jacobian matrix of the vector field $f$ has the form
$$
f'(p_{0})=\left(\begin{array}{cc} {-ab} & {0} \\ {0} & {-cd}
\end{array}\right)=\left(\begin{array}{cc} -2{\sqrt{3}} & {0} \\ {0}
& -2{\sqrt{3}}
\end{array}\right),$$  $$f'(p_{1})=\left(\begin{array}{cc} {b(a-b)} & {0} \\ {0}
& {d(c-d)}
\end{array}\right)=\left(\begin{array}{cc} {2\sqrt{3}-3} & {0} \\ {0}
& {2\sqrt{3}-2}
\end{array}\right),
$$
$$
f'(p_{2})=\left(\begin{array}{cc} {a(b-a)} & {0} \\
{0} & {c(d-c)}
\end{array}\right)=\left(\begin{array}{cc} {2\sqrt{3}-4} & {0} \\
{0} & {2\sqrt{3}-6}
\end{array}\right),
$$
$$
f'(p_{3})=\left(\begin{array}{cc} {b(a-b)} & {0} \\ {0} &
{c(d-c)}
\end{array}\right)=\left(\begin{array}{cc} {2\sqrt{3}-3} & {0} \\ {0} &
{2\sqrt{3}-6}
\end{array}\right).
$$
We see that $\delta(p_{i})=i$, $0\leq i\leq 3$, and hence this
system does not admit a normally hyperbolic inertial manifold in
state space $X$ by Lemma~3.3. $\Box$

REMARK 3.6. Thus, we have separated the classes of problems
admitting inertial manifolds and normally hyperbolic inertial
manifolds for 3D two-component systems of chemical kinetics
equations with a cubic nonlinearity. In particular, we have
obtained an inertial manifold that is not normally hyperbolic.

\bigskip
\centerline{\bf 4.  Absolutely Normally Hyperbolic Inertial
Manifolds}
\bigskip

Under assumptions (i)--(iii), the same spectrum sparseness
condition~(2.2) is sufficient for the existence of an absolutely
normally hyperbolic inertial manifold ${\EuScript M}\subset
X^{\theta}$ for an arbitrary nonlinear part $F(u)$ of Eq.~(1.1)
(see~[17, Theorem~5.6] and [18, Theorem~81.4]).\footnote{Such
manifolds are called normally hyperbolic in~[17,~18].}

Consider scalar homogeneous equations of the form
$$
 \partial _{t} u \, =\nu\Delta u+f(u),\; \; \nu>0, \eqno(4.1)
$$
in a bounded Lipschitz domain $\Omega\subset \mathbb{R}^{m}$,
$m\leq3$, with the Neumann condition $(\rm{N})$ or the periodicity
condition $(\rm{P})$ on $\partial \Omega$ and with a function $f\in
C^{3}(\mathbb{R},\mathbb{R})$ satisfying the sign condition. Let
$\sigma(-\Delta)=\{0=\lambda_{1}\leq \lambda_{2}\leq\dotsm\}$.
Being a special case of~(2.3), the dissipative equation~(4.1) can
be represented in the form~(1.1) with $X=L^{2}(\Omega)$ and with
the nonlinearity exponent $\theta=0$.

LEMMA 4.1 ([15]). \textit{Let $\lambda_{n+1}-\lambda_{n}\leq K$,
$n\geq 1$, and let $f'(p_{0})-f'(p_{1})=a>0$ for some
$p_{0},p_{1}\in \mathbb{R}$ such that $f(p_{0})=f(p_{1})=0.$ Then
problems $(4.1)_{\rm{N}}$, and~$(4.1)_{\rm{P}}$ do not have a
normally hyperbolic inertial manifold ${\EuScript M}\subset X$ for
$\nu < a/K$}.

A simple proof is based on Lemma~3.2.

COROLLARY 4.2. \textit{If $\lambda_{n+1}-\lambda_{n}\leq K$, $n\geq
1$, and $f(u)=u-u^{3}$, then Eq.~$(4.1)$ with the boundary
condition $(\rm{N})$ or $(\rm{P})$ does not have an absolutely
normally hyperbolic inertial manifold ${\EuScript M}\subset X$ for
$\nu < 3/K$}.

In the case of $\Omega=(0,\pi)^{3}$, the spectrum of the operator
$-\Delta$ with the Neumann condition or the periodicity condition
on $\partial \Omega$ consists of eigenvalues of the form
$\lambda_{n}=l_{1}^{2}+l_{2}^{2}+l_{3}^{2}$, $l_{j}\in \mathbb{Z}$;
here one always has $\lambda_{n+1}-\lambda_{n}\leq3$ by the Gauss
theorem [3], and hence one can take $K=3$ in Corollary~4.2.

COROLLARY 4.3. \textit{Equation $(4.1)$ with $f(u)=u-u^{3}$ and
with one of the boundary conditions $(\rm{N})$ and $(\rm{P})$ in
the cube $\Omega=(0,\pi)^{3}$ does not have an absolutely normally
hyperbolic inertial manifold ${\EuScript M}\subset X$ for $\nu<1$}.

We see that an absolutely normally hyperbolic inertial manifold may
fail to exist even for very simple semilinear parabolic equations.

\bigskip
\centerline{\bf 5. Conclusion}
\bigskip

As was already mentioned, the technique in~[6] permits one to
derive the existence of a normally hyperbolic inertial manifold in
an appropriate state space for Eqs.~(2.3) and~(4.1) from the
principle of spatial averaging for the Laplace operator. Since this
principle holds for the 3D cube, a careful solution of this problem
will (in view of Corollary~4.3) permit separating the classes of
problems admitting a normally hyperbolic inertial manifold and an
absolutely normally hyperbolic inertial manifold for 3D scalar
chemical kinetics equations.

There is a suspicion that, for an appropriate choice of the phase
space and the family of admissible nonlinearities, the validity of
the principle of spatial averaging and the sparseness of the
spectrum of the Laplace operator in the scalar reaction--diffusion
equations are \textit{necessary and sufficient} for the existence
of a normally hyperbolic inertial manifold and an absolutely
normally hyperbolic inertial manifold, respectively. Needless to
say, we speak of the existence of such manifolds for every
nonlinearity in a given family.

CONJECTURE 5.1. \textit{The following properties are equivalent for
equations of the form~$(2.3)$ in a bounded Lipschitz domain
$\Omega\subset \mathbb{R}^{m}$, $m\leq3$, with the boundary
conditions $(\rm{D})$, $(\rm{N})$, or $(\rm{P})$ }:

(a) \textit{The validity of the principle of spatial averaging for
the Laplace operator $\Delta_{\Omega}$}.

(b) \textit{The existence of a normally hyperbolic inertial
manifold in an appropriate state space for an arbitrary
``admissible'' function $f$ and an arbitrary diffusion
coefficient~$\nu$}.

The implication (a) $\Rightarrow$ (b) can be derived by the
technique in [6, 20] under the assumption of sufficient smoothness
of the operator $u\rightarrow f(x,u)$ in the corresponding
functional space. The main problem is to establish the converse
implication for the right choice of the family of admissible
functions $f$.

CONJECTURE 5.2. \textit{The following properties are equivalent for
equations of the form~$(4.1)$ in a bounded Lipschitz domain
$\Omega\subset \mathbb{R}^{m}$, $m\leq3,$ with the boundary
conditions $((\rm{N})$ or $(\rm{P})$}:

(a) \textit{The sparseness of the spectrum of the Laplace
operator~$\Delta_{\Omega}$}.

(b) \textit{The existence of an absolutely normally hyperbolic
inertial manifold in an appropriate state space for an arbitrary
``admissible'' function $f$ and an arbitrary diffusion
coefficient~$\nu$}.

The implication (b) $\Rightarrow$ (a) follows from Corollary~4.2,
provided that admissible functions include cubic polynomials. The
converse can be obtained by the technique in~[18, Section~5] if
one starts from a ``smoother'' main space of the parabolic
equation, say, by setting $X=H^{s}(\Omega)$ for some $s\geq1$.

\bigskip
\medskip \centerline{\bf REFERENCES}
\medskip

\noindent 1. S.-N. Chow, K. Lu, and G. R. Sell, ``Smoothness of
inertial manifolds, ''\textit{J. Math. Anal. Appl.}, \textbf{169},
no. 1, (1992), 283--312.

\noindent 2. A. Eden, V. Kalantarov, and S. Zelik,
``Counterexamples to the regularity of {M}a{\~n}{\'e} projections
in the attractors theory,'' \textit{Russian Math. Surveys},
\textbf{68}, no. 2, (2013), 199--226.

\noindent 3. G. H. Hardy and E. M. Wright, \textit{An Introduction
to the Theory of Numbers}, 5th ed., Oxford Press, Oxford, 1979.

\noindent 4. D. Henry, \textit{Geometric Theory of Semilinear
Parabolic Equations}, Lecture Notes in Math., \textbf{840},
Springer, New York, 1981.

\noindent 5. M. Hirsch, G. Pugh, and M. Shub, \textit{Invariant
Manifolds}, Lecture Notes in Mathematics, \textbf{583}, Springer,
New York, 1977.

\noindent 6. A. Kostianko and S. Zelik, "Inertial manifolds for 3D
Cahn--Hilliard equations with periodic boundary conditions,''
\textit{Comm. Pure Appl. Anal.}, \textbf{14}, no. 5, (2015),
2069--2094.

\noindent 7. A. Kostianko, ``Inertial manifolds for the 3D modified
Leray-$\alpha$ model with periodic boundary conditions,''
arXiv:1510.08936v1.

\noindent 8. H. Kwean, ``An extension of the principle of spatial
averaging for inertial manifolds,''  \textit{J. Austral. Math.
Soc., Ser. A}, \textbf{66}, no. 1, (1999), 125--142.

\noindent 9. H. Kwean, ``A geometric criterion for the weaker
principle of spatial averaging,'' \textit{Comm. Korean Math. Soc.},
\textbf{14}, no. 2, (1999), 337--352.

\noindent 10. J. Mallet-Paret and G. R. Sell, ``The principle of
spatial averaging and inertial manifolds for reaction-diffusion
equations,'' Lecture Notes in Math., \textbf{1248}, Springer, New
York, 1987, pp. 94--107.

\noindent 11. J. Mallet-Paret and G. R. Sell, ``Inertial manifolds
for reaction diffusion equations in higher space dimensions,''
\textit{J. Amer. Math. Soc.}, \textbf{1}, no. 4, (1988), 805--866.

\noindent 12. J. Mallet-Paret, G. R. Sell, and Z. Shao,
``Obstructions to the existence of normally hyperbolic inertial
manifolds,''  \textit{Indiana Univ. Math. J.}, \textbf{42}, no. 3,
(1993), 1027--1055.

\noindent 13. V. A. Pliss and G. R. Sell, ``Perturbations of
normally hyperbolic manifolds with applications to the
Navier--Stokes equations,'' \textit{J. Diff. Equat.}, \textbf{169},
no. 2, (2001), 396--492.

\noindent 14. J. C. Robinson, \textit{Infinite-Dimensional
Dynamical Systems}, Cambridge Texts in Applied Mathematics,
Cambridge University Press, Cambridge, 2001.

\noindent 15. A. V. Romanov, ``Three counterexamples in the theory
of inertial manifolds,'' \textit{Math. Notes}, \textbf{68}, no.
3--4, 2000,  378--385.

\noindent 16. A. V. Romanov, ``Parabolic equation with nonlocal
diffusion without a smooth inertial manifold,'' \textit{Math.
Notes}, \textbf{96}, no. 4, (2014), 548--555.

\noindent 17. R. Rosa and R. Temam, ``Inertial manifolds and normal
hyperbolisity,'' \textit{ACTA Applicandae Mathematicae},
\textbf{45} (1996), 1--50.

\noindent 18. G. R. Sell and Y. You, \textit{Dynamics of
Evolutionary Equations}, Appl. Math. Sci., \textbf{143}, Springer,
New York, 2002.

\noindent 19. R. Temam, \textit{Infinite-Dimensional Dynamical
Systems in Mechanics and Physics}, 2d ed., Appl. Math. Sci.,
\textbf{68}, Springer, New York, 1997.

\noindent 20. S. Zelik, ``Inertial manifolds and finite-dimensional
reduction for dissipative PDEs,'' \textit{Proc. Roy. Soc.
Edinburgh, Ser. A}, \textbf{144}, no. 6, (2014), 1245--1327.

\end{document}